\documentclass[12pt]{article}

\usepackage{amsfonts}
\usepackage{amsmath}
\usepackage{tocbibind}
\newtheorem{thm}{Theorem}[section]
\newtheorem{rem}{Remark}[section]

\newcommand{\dC}{\mathbb{C}}

\newcommand{\dE}{\mathbb{E}}
\newcommand{\dP}{\mathbb{P}}

\newcommand{\cB}{\mathcal{B}}

\newcommand{\cN}{\mathcal{N}}

\newcommand{\cF}{\mathcal{F}}

\newcommand{\rI}{\mathrm{I}}
\newcommand{\veps}{\varepsilon}

\newcommand{\ind}{\mbox{1}\kern-.25em \mbox{I}}
\font\calcal=cmsy10 scaled\magstep1
\def\build#1_#2^#3{\mathrel{\mathop{\kern 0pt#1}\limits_{#2}^{#3}}}
\def\liml{\build{\longrightarrow}_{}^{{\mbox{\calcal L}}}}
\def\limp{\build{\longrightarrow}_{}^{{\mbox{\calcal P}}}}
\def\videbox{\mathbin{\vbox{\hrule\hbox{\vrule height1ex \kern.5em
\vrule height1ex}\hrule}}}
\newtheorem{lem}{Lemma}[section]

\newtheorem{cor}[lem]{Corollary}

\def\demend{\hfill $\videbox$\\}

\usepackage{verbatim}

\renewcommand{\theequation}{\thesection.\arabic{equation}}

\begin{document}

\title{New insights on the minimal random walk}

\author{Bernard Bercu$^1$ and V\'ictor Hugo V\'azquez Guevara$^2$}

\date{ }
	\maketitle
$^1$Universit\'e de Bordeaux, Institut de Math\'ematiques de Bordeaux, 351 cours de la 
lib\'eration, 33405 Talence, France\\

$^2$Benem\'erita Universidad Aut\'onoma de Puebla, Facultad de Ciencias F\'isico
Matem\'aticas, Avenida San Claudio y R\'io Verde, 72570 Puebla, M\'exico.\\

$^1$bernard.bercu@math.u-bordeaux.fr, $^2$vvazquez@fcfm.buap.mx
\vspace{10pt}

\begin{abstract}
The aim of this paper is to deepen the analysis of the asymptotic behavior of
the so-called minimal random walk (MRW) using a new martingale approach. 
The MRW is a discrete-time random walk with infinite memory that has
three regimes depending on the location of its two parameters.
In the diffusive and critical regimes, we establish new results on the almost sure asymptotic
behavior of the MRW, such as the quadratic strong law and the law of the
iterated logarithm. In the superdiffusive regime, we prove the almost sure convergence of the 
MRW, properly normalized, to a nondegenerate random variable. Moreover, we show that the fluctuation
of the MRW around its limiting random variable is still Gaussian.
\end{abstract}

%
\vspace{2pc}

\section{Introduction}

The minimal random walk (MRW) was first proposed by Harbola, Kumar and Lindenberg 
\cite{Kumar2014} in $2014$.  It can be seen as a variant of the famous elephant random walk (ERW)
introduced in the early 2000s by Sch\"utz and Trimper \cite{Schutz}, in order to investigate
how long-range memory affects the behavior of the random walk, see also \cite{Baur2016, Bercu2018, Coletti2017a, Coletti2017b,Vazquez2019}. 
The movements of the walker in
the MRW are quite simple with only two possibilities, either a forward step or a
resting step, depending on the choice of two parameters $p$ and $q$ in $[0,1]$.
\vspace{1ex} \\
The MRW  is defined as follows. The walker is located at the origin at time zero, $S_0=0$. For the first step, $S_1=X_1$ where $X_1$
has the Bernoulli $\cB(s)$ distribution which means that the walker
goes to the right at point 1 with probability $s$ or stays at the origin with probability 
$1-s$ for some $s$ in $[0,1]$. Afterwards, at time $n+1\geq 2$, an integer $k$ is chosen
uniformly at random among the previous times $1,\ldots,n$ and the step $X_{n+1}$ is determined
stochastically by
\begin{equation}
   X_{n+1} = \left \{ \begin{array}{ccc}
    \alpha_{n+1} &\text{ if } & X_k=1, \vspace{1ex}\\
    \beta_{n+1} &\text{ if } & X_k=0,
   \end{array} \right.
   \nonumber
\end{equation}
where $\alpha_{n+1}$ and $\beta_{n+1}$ are two independent discrete random variables
with Bernoulli $\cB(p)$, and $\cB(q)$ distributions, respectively.  In other words,
\begin{equation}
\label{STEPS}
   X_{n+1} = \alpha_{n+1} X_{U_n}+\beta_{n+1}(1-X_{U_n})
\end{equation}
where $U_n$ is a discrete uniform random variable on $\bigl\{1,\ldots,n\bigr\}$
and $\alpha_{n+1}$, $\beta_{n+1}$ and $U_n$ are mutually independent. Then, the position of the MRW is
given by
\begin{equation}
\label{POSMRW}
S_{n+1}=S_{n}+X_{n+1}.
\end{equation}
The fundamental parameter of the MRW is defined by the difference
\begin{equation}
\label{DEFA}
a=p-q.
\end{equation}
Throughout the paper, we assume that $a<1$ inasmuch as $a =1$ only appears in the
trivial case where $p =1$ and $q= 0$ which means that for all $n \geq 1$, $X_n=X_1$.
The MRW is said to be diffusive if $a<1/2$, critical if $a=1/2$ and supercritical if $a>1/2$. It has been recently shown by
Coletti, Gava and Lima \cite{Coletti2019} that whatever the value of the parameter $a$ in $[-1,1[$,
\begin{equation}
\label{ASCVGSN}
\lim_{n \rightarrow \infty} \frac{S_n}{n}=\frac{q}{1-a}
\hspace{1cm}\text{a.s.}
\end{equation}
Moreover, it has also been proven in \cite{Coletti2019} that in the diffusive regime $a< 1/2$,
\begin{equation}
\label{AND}
\sqrt{n}\Bigl(\frac{S_n}{n}-\frac{q}{1-a}\Bigr) \liml N\Bigl(0,\frac{\sigma^2}{1-2a}\Bigr)
\end{equation}
while in the critical regime $a= 1/2$,
\begin{equation}
\label{ANC}
\sqrt{\frac{n}{\log n}}\Bigl(\frac{S_n}{n}-2q \Bigr)\liml N\bigl(0,\sigma^2\bigr)
\end{equation}
where the asymptotic variance
$$
\sigma^2=\frac{q(1-p)}{(1-a)^2}.
$$
Laws of iterated logarithm were also established for both diffusive and critical regimes. Furthermore,
in the superdiffusive regime $a>1/2$ and in the special situation where $q=0$ and $p>1/2$, it has been shown
in \cite{Coletti2019} that
\begin{equation}
\label{ASCVGS}
 \lim_{n \rightarrow \infty} \frac{S_n}{n^{p}}=L \hspace{1cm} \text{a.s.}
\end{equation}
where $L$ is a non-degenerate random variable which is non-Gaussian. 
Finally, by the calculation of all factorial moments of $S_n$, it
has been recently proven by Miyazaki and Takei \cite{Miyazaki2020} that in the superdiffusive regime $a>1/2$ 
with $q=0$ and $p>1/2$, $L$ has a Mittag-Leffler distribution with parameter $p$. It implies in particular that
$\dP(L>0)=1$.
\vspace{-1ex} \\
The aim of this paper is to deepen the analysis of \cite{Coletti2019} and \cite{Miyazaki2020} in several directions. 
We shall make use of an alternative martingale approach, similar to the one used 
by Bercu in \cite{Bercu2018} for the ERW,
which allows us to carry out the asymptotic analysis of the MRW in a more natural way.
On the one hand, we shall prove functional central limit theorems both diffusive and critical regimes, which in turn
imply corresponding central limit theorems \cite{Coletti2019}. 
Beside, we shall also establish new quadratic strong laws which are
really useful in statistical applications \cite{Bercu2021b}.
On the other hand, we shall extend the results of \cite{Coletti2019} and \cite{Miyazaki2020} in the superdiffusive regime 
$a>1/2$ without assuming that $q=0$.  In particular, as it was already done
for the ERW \cite{Kubota2019}, we shall show that the fluctuation of the MRW around 
its limiting random variable is still Gaussian.
\vspace{1ex} \\
The paper is organized as follows. Section \ref{S-MR} is devoted to the main results of the paper.
Our first contribution is to establish functional central limit theorems for the MRW 
in the diffusive and critical regimes. Our second contribution is also to provide 
a functional central limit theorem and to investigate the fluctuation of the MRW around 
its limiting random variable in the superdiffusive regime. 
Four Appendices are included. Our martingale approach,
slightly different from that of \cite{Coletti2019, Miyazaki2020}, is described in Appendix A.
Appendices B, C and D contain the proofs of the results in the diffusive, critical and superdiffusive regimes, respectively.

\section{Main Results}
\label{S-MR}
\setcounter{equation}{0}

This section is devoted to the main results on the asymptotic behavior of the MRW.

\subsection{The diffusive regime}

Our first results concern to the asymptotic behavior of the MRW in the diffusive regime
where $a<1/2$.
The strong law of large numbers and the law of iterated logarithm were previously established in 
Theorems 1 and 3 of \cite{Coletti2019}, respectively. The quadratic strong law is new. In all the sequel, we will make use of the asymptotic variance
\begin{equation}
\label{VAR-DR}
\sigma^2=\frac{q}{1-a}\Big(1-\frac{q}{1-a} \Big)=\frac{q(1-p)}{(1-a)^2}.
\end{equation}

\begin{thm}
\label{T-ASP-DR}
We have the almost sure convergence
\begin{equation}
\label{SLLN-DR}
\lim_{n \rightarrow \infty} \frac{S_n}{n}=\frac{q}{1-a} \hspace{1cm} \text{a.s.}
\end{equation}
In addition, we also have the law of iterated logarithm
\begin{eqnarray}
\hspace{-2.2cm}
\limsup_{n \rightarrow \infty}  \Big(\frac{n}{2 \log \log n}\Big)^{1/2} 
\Big(\frac{S_n}{n}-\frac{q}{1-a}\Big) 
&=& -\liminf_{n \rightarrow \infty}  \Big(\frac{n}{2 \log \log n}\Big)^{1/2} 
\Big(\frac{S_n}{n}-\frac{q}{1-a}\Big) \nonumber \\
&=& \frac{\sigma}{\sqrt{1-2a}} \hspace{1cm}\text{a.s.}
\label{LIL-DR}
\end{eqnarray}
In particular, 
\begin{equation}
\label{LILS-DR}
\limsup_{n \rightarrow \infty}  \Big(\frac{n}{2 \log \log n}\Big) 
\Big(\frac{S_n}{n}-\frac{q}{1-a}\Big)^2=\frac{\sigma^2}{1-2a}\hspace{1cm}\text{a.s.}
\end{equation}
Moreover, we have the quadratic strong law 
\begin{equation}
\label{QSL-DR}
\lim_{n \rightarrow \infty} \frac{1}{\log n}\sum_{k=1}^n \Big(\frac{S_k}{k}-\frac{q}{1-a}\Big)^2 = \frac{\sigma^2}{1-2a}\hspace{1cm}\text{a.s.}
\end{equation}
\end{thm}

\noindent
Hereafter, we focus our attention on the distributional convergence of the MRW.
Denote by $D([0,\infty[)$ the Skorokhod space of right-continuous functions with left-hand limits.
Our functional central limit theorem, which extends Theorem 2 in \cite{Coletti2019}, is as follows.

\begin{thm}
\label{T-FCLT-DR}
We have the distributional convergence in $D([0,\infty[)$,
\begin{equation}
\label{FCLT-DR}
\left( \sqrt{n}\Big(\frac{S_{\lfloor nt \rfloor}}{\lfloor nt \rfloor}-\frac{q}{1-a}\Big), t \geq 0\right) \Longrightarrow \big( W_t, t \geq 0 \big)
\end{equation}
where $\big( W_t, t \geq 0 \big)$ is a real-valued centered Gaussian process starting at the origin with covariance given, for all $0<s \leq t$, by
$$
\dE[W_s W_t]= \frac{\sigma^2}{(1-2a)t} \Big(\frac{t}{s}\Bigr)^a.
$$
In particular, we have the asymptotic normality
\begin{equation}
\label{CLT-DR}
\sqrt{n}\Big(\frac{S_n}{n}-\frac{q}{1-a}\Big) \liml N\Big(0,\frac{\sigma^2}{1-2a}\Big).
\end{equation}
\end{thm}

\noindent
It is also interesting to investigate the asymptotic behavior of the center of mass of the
MRW defined by

$$
G_n=\frac{1}{n}\sum_{k=1}^n S_k.
$$

\noindent
Very recent results on the center of mass of the ERW can be found in \cite{Bercu2021}. 
The strong law of large numbers for $(G_n)$ follows from \eqref{SLLN-DR} while the asymptotic
normality is a direct application of the distributional convergence \eqref{FCLT-DR}.

\begin{cor}
\label{C-CM-DR}
We have the almost sure convergence
\begin{equation}
\label{SLLN-CM-DR}
\lim_{n \rightarrow \infty} \frac{G_n}{n}=\frac{q}{2(1-a)} \hspace{1cm} \text{a.s.}
\end{equation}
Moreover, we have the asymptotic normality
\begin{equation}
\label{CLT-CM-DR}
\sqrt{n}\Big(\frac{G_n}{n}-\frac{q}{2(1-a)}\Big) \liml N\Big(0,\frac{2\sigma^2}{3(1-2a)(2-a)}\Big).
\end{equation}
\end{cor}

\subsection{The critical regime}

We now study the asymptotic behavior of the MRW in the critical regime
where $a=1/2$. Once again, the strong law of large numbers and the law of iterated logarithm were previously established in Theorems 1 and 3 of \cite{Coletti2019}, while the quadratic strong law is new.

\begin{thm}
\label{T-ASP-CR}
We have the almost sure convergence
\begin{equation}
\label{SLLN-CR} 
\lim_{n \rightarrow} \frac{S_n}{n}=2q  \hspace{1cm} \text{a.s.}
\end{equation}
In addition, we also have the law of iterated logarithm
\begin{eqnarray}
& &\limsup_{n \rightarrow \infty}  \Big(\frac{n}{2 \log n \log \log \log n}\Big)^{1/2} 
\Big(\frac{S_n}{n}-2q\Big) \nonumber \\
& = & -\liminf_{n \rightarrow \infty}  \Big(\frac{n}{2 \log n \log \log \log n}\Big)^{1/2} 
\Big(\frac{S_n}{n}-2q\Big) 
\nonumber \\
&=&  \sqrt{4q(1-p)} \hspace{1cm}\text{a.s.}
\label{LIL-CR}
\end{eqnarray}
In particular, 
\begin{equation}
\label{LILS-CR}
\limsup_{n \rightarrow \infty}  \Big(\frac{n}{2 \log n \log \log \log n}\Big)
\Big(\frac{S_n}{n}-2q\Big)^2=4q(1-p)\hspace{1cm}\text{a.s.}
\end{equation}
Moreover, we have the quadratic strong law 
\begin{equation}
\label{QSL-CR}
\lim_{n \rightarrow \infty} \frac{1}{\log \log n}\sum_{k=1}^n 
\Big( \frac{1}{\log k} \Big)^2\Big(\frac{S_k}{k}-2q\Big)^2 = 4q(1-p)\hspace{1cm}\text{a.s.}
\end{equation}
\end{thm}

Our next result concerns the functional central limit theorem in the critical regime.

\begin{thm}
\label{T-FCLT-CR}
We have the distributional convergence in $D([0,\infty[)$,
\begin{equation}
\label{FCLT-CR}
\left( \sqrt{\frac{n^t}{\log n}}\Big(\frac{S_{\lfloor n^t \rfloor}}{\lfloor n^t \rfloor}-2q\Big), t \geq 0\right) \Longrightarrow \big( 2\sqrt{q(1-p)} B_t, t \geq 0 \big)
\end{equation}
where $\big( B_t, t \geq 0 \big)$ is a standard Brownian motion.
In particular, we have the asymptotic normality
\begin{equation}
\label{CLT-CR}
\sqrt{\frac{n}{\log n}}\Big(\frac{S_n}{n}-2q\Big) \liml N\big(0,4q(1-p)\big).
\end{equation}
\end{thm}

\noindent
The asymptotic behavior of the center of mass in the critical regime is as follows.

\begin{cor}
\label{C-CM-CR}
We have the almost sure convergence
\begin{equation}
\label{SLLN-CM-CR}
\lim_{n \rightarrow \infty} \frac{G_n}{n}=q \hspace{1cm} \text{a.s.}
\end{equation}
Moreover, we have the asymptotic normality
\begin{equation}
\label{CLT-CM-CR}
\sqrt{\frac{n}{\log n}}\Big(\frac{G_n}{n}-q\Big) \liml N\Big(0,\frac{16q(1-p)}{9}\Big).
\end{equation}
\end{cor}

\subsection{The superdiffusive regime}

The superdiffusive regime is more difficult to handle as it 
requires more technical considerations. 
It has already been treated in Theorem 4 of \cite{Coletti2019} but only in the special case where 
$q=0$ and $p>1/2$. Our aim is now to extend the previous results of 
\cite{Coletti2019, Miyazaki2020} to the general case where $q\geq 0$ and $a>1/2$.

\begin{thm}
\label{T-ASP-SR}
We have the almost sure convergence
\begin{equation}
\label{FSLLN-SR}
\left( n^{1-a}\Big(\frac{S_{\lfloor nt \rfloor}}{\lfloor nt \rfloor}-\frac{q}{1-a}\Big), t > 0\right) \longrightarrow \Big( \frac{1}{t^{1-a}}L, t > 0 \Big)
\end{equation}
where $L$ is a non-degenerated random variable. In particular
\begin{equation}
\label{SLLN-SR}
\lim_{n \rightarrow \infty} n^{1-a}\Big(\frac{S_n}{n}-\frac{q}{1-a}\Big) = L \hspace{1cm}\text{a.s.}
\end{equation}
Moreover, this convergence also holds in $L^2$, 
\begin{equation} 
\label{SDL2-SR}
\lim_{n\rightarrow \infty} \dE \Big[ n^{1-a}\Big(\frac{S_n}{n}-\frac{q}{1-a}-L\Big)^2 \Big]=0.
\end{equation}
\end{thm}

\begin{rem}
Since $a=p-q$, $1-a=1-p+q$. Consequently, \eqref{SLLN-SR} implies that
\begin{equation} 
\lim_{n \rightarrow} n^{1-p+q}\Big(\frac{S_n}{n}-\frac{q}{1-p+q}\Big) = L \hspace{1cm}\text{a.s.}
\nonumber
\end{equation}
One can observe that in the 
special case where $q=0$, we find again 
convergence \eqref{ASCVGS}.
\end{rem}

\begin{thm}
\label{T-MOML}
The first two moments of the limiting random variable $L$ are given by
\begin{equation}
\dE[L]=\frac{s+\nu}{\Gamma(p-q+1)} \hspace{1cm} \text{and}
\hspace{1cm}
\dE[L^2]=\frac{s+\tau}{\Gamma(2(p-q)+1)},
\end{equation}
where $\nu=-q/(1-a)$ and $\tau$ is given by
\begin{equation}
\label{deftau}
\tau=s-\frac{4qs}{1-a}+\frac{2q}{2a-1}+
\frac{4q^2(3a-2)}{(1-a)^2 (2a-1)}.
\end{equation}
\end{thm}

\noindent
We now focus our attention on the fluctuation of the MRW around its limiting random variable $L$, 
in the spirit of the original work of Kubota and Takei \cite{Kubota2019}. 

\begin{thm}
\label{T-AN-SR}
We have the asymptotic normality
\begin{equation}
\label{AN-SR}
\sqrt{n^{2a-1}}\left( n^{1-a}\Big(\frac{S_n}{n}-\frac{q}{1-a}\Big)-L \right) 
\liml \cN\Big(0,\frac{\sigma^2}{2a-1}\Big).
\end{equation}
\end{thm}

\begin{rem}
One can observe that the fluctuation around $L$ is still Gaussian. 
Moreover, one can notice that the asymptotic variance coincides with the one
obtained for the central limit theorem \eqref{CLT-DR} in the diffusive regime $a<1/2$, up to a sign factor.
Finally, the asymptotic normality \eqref{AN-SR} was stated without proof in \cite{Miyazaki2020} 
where only the special case $q=0$ was taken into consideration.
\end{rem}

\noindent
Our last result concerns the center of mass in the superdiffusive regime.

\begin{cor}
\label{C-CM-SR}
We have the almost sure convergence
\begin{equation}
\label{SLLN-CM-SR}
\lim_{n \rightarrow \infty} n^{1-a}\Big(\frac{G_n}{n}-\frac{q}{2(1-a)}\Big)=\frac{L}{1+a} \hspace{1cm} \text{a.s.}
\end{equation}
\end{cor}

\section*{Appendix A. Our martingale approach}
\renewcommand{\thesection}{\Alph{section}}
\renewcommand{\theequation}{\thesection.\arabic{equation}}
\setcounter{section}{1}
\setcounter{equation}{0}
\vspace{-1ex}

It follows from \eqref{STEPS} that for all $n\geq 1$,
\begin{equation}
\dE[X_{n+1}|\cF_n]=\dE[\alpha_{n+1}]\dE[X_{U_n}|\cF_n]+\dE[\beta_{n+1}] (1-\dE[X_{U_n}|\cF_n] ) \hspace{1cm}\text{a.s.}
 \nonumber 
\end{equation}
where $\cF_n$ is the natural $\sigma$-algebra, $\mathcal{F}_n=\sigma\left(X_1,\ldots,X_n\right)$. Consequently, as
$U_n$ is a discrete uniform random variable on $\{1,\ldots,n\}$, we obtain that
\begin{equation}
\label{e1}
\dE[X_{n+1}|\cF_n]=p \frac{S_n}{n}+q\Big(1 - \frac{S_n}{n}\Big)=q+a \frac{S_n}{n} \hspace{1cm}\text{a.s.}
\end{equation}
Hence, \eqref{POSMRW} together with \eqref{e1} imply that almost surely
\begin{equation}
\label{expos}
\dE[S_{n+1}|\cF_n]=q+\gamma_n S_n \hspace{1cm}\text{where}\hspace{1cm}
\gamma_n=1+\frac{a}{n}.
\end{equation}
Let $(a_n)$ be the deterministic sequence given by $a_1=1$ and, for $n\geq 2$,
\begin{equation}
\label{an}
a_n=\prod_{k=1}^{n-1}\gamma_k^{-1}=\frac{\Gamma(n)\Gamma(a+1)}{\Gamma(n+a)},
\end{equation}%
where $\Gamma$ stands for the Euler Gamma function.
In order to define the martingale $(M_n)$ that will lead us to the asymptotic analysis of the MRW, 
let us introduce the sequence $(A_n)$ given by $A_0=0$ and, for $n\geq 1$,
\begin{equation}
\label{bigAn}
A_n=\sum_{k=1}^n a_n.
\end{equation}
Denote $M_0=0$ and, for $n\geq 1$,
\begin{equation} 
\label{martingale}
M_n=a_n S_n-q A_n.
\end{equation}
We clearly deduce from \eqref{expos}, \eqref{an} and \eqref{bigAn} that almost surely
$$ 
\dE[M_{n+1}|\mathcal{F}_n]=a_{n+1}(q+\gamma_n S_n)-q A_{n+1} = a_{n+1}\gamma_n S_n-q A_n =a_n S_n-q A_n=M_n
$$
It means that the sequence $(M_n)$ is a discrete-time martingale such that, for all $n\geq 1$,
$\dE[M_n]=\dE[M_1]=s-q$. One can observe that our martingale $(M_n)$ is
slightly different from that of \cite{Coletti2019, Miyazaki2020}. Our aim is to investigate the asymptotic behavior of $(M_n)$ in order to deduce the asymptotic behavior of the position $(S_n)$ of the MRW.
The martingale $(M_n)$ can be rewritten in the additive form
\begin{equation} \label{martingale2}
M_n=\sum_{k=1}^n a_k \varepsilon_k
\end{equation}
where, for $n\geq 1$, the martingale increment $\varepsilon_n=S_n-\dE[S_n|\mathcal{F}_{n-1}]=S_{n}-(q+\gamma_{n-1} S_{n-1})$. 
The predictable quadratic variation associated with $(M_n)$ is given by
$\langle M \rangle_0=0$ and, for all $n \geq 1$,
\begin{equation} 
\label{IPM}
\langle M \rangle_n = \sum_{k=1}^n a_k^2\dE[\varepsilon_k^2| \mathcal{F}_{k-1}].
\end{equation}
By recalling that $X_{n+1}=X_{n+1}^2$, we immediately have from \eqref{e1} that
\begin{equation}
\label{e2}
\dE[X_{n+1}^2|\cF_n]=q+a \frac{S_n}{n} \hspace{1cm}\text{a.s.}
\end{equation}
Hence, we deduce from \eqref{POSMRW} and \eqref{e2} that
\begin{eqnarray}
\label{Snsqq}
\mathbb{E}[S_{n+1}^2|\mathcal{F}_{n}] &=& S_n^2+2S_n\Big( q+a\frac{S_n}{n} \Big)+\Big( q+a\frac{S_n}{n} \Big) \hspace{1cm} \text{a.s.} \notag  \\
&=&S_n^2+(1+2S_n)\Big( q+a\frac{S_n}{n} \Big) \hspace{1cm} \text{a.s.}
\end{eqnarray}
Therefore, as $\mathbb{E}[\varepsilon^2_{n+1}|\mathcal{F}_n]= \mathbb{E}[S_{n+1}^2|\mathcal{F}_{n}] - (q+\gamma_{n} S_{n})^2$, 
we obtain from \eqref{Snsqq} that 
\begin{eqnarray} 
\mathbb{E}[\varepsilon^2_{n+1}|\mathcal{F}_n]
&=&
S_n^2+(1+2S_n)\Big( q+a\frac{S_n}{n} \Big)-\Big(q+S_n+a\frac{S_n}{n}\Big)^2  \hspace{1cm} \text{a.s.}\notag \\
&=&
\Big(q+a\frac{S_n}{n} \Big)-\Big(q+a\frac{S_n}{n} \Big)^2 \hspace{1cm} \text{a.s.}
\label{M2EPS}
\end{eqnarray}
Equation \eqref{M2EPS} clearly leads to
\begin{equation}
\label{SUPEPS2}
\sup_{n\geq 0} \mathbb{E}\left[\varepsilon^2_{n+1}|\mathcal{F}_n\right]\leq \frac{1}{4} \hspace{1cm} \text{a.s.}
\end{equation}
On the same direction, we also have
\begin{equation}
\dE[X_{n+1}^3|\cF_n]=\dE[X_{n+1}^4|\cF_n]=q+a \frac{S_n}{n} \hspace{1cm}\text{a.s.} \notag
\end{equation}
which implies that
\begin{eqnarray*}
\mathbb{E}[S_{n+1}^3|\mathcal{F}_n] &=& S_n^3+\big(1+3S_n+3S_n^2\big)\Big(q+a\frac{S_n}{n}\Big)\hspace{1cm}\text{a.s.} \\
\mathbb{E}[S_{n+1}^4|\mathcal{F}_n] &=& S_n^4+\big(1+4S_n+6S_n^2+4S_n^3\big)\Big(q+a\frac{S_n}{n}\Big)\hspace{1cm}\text{a.s.}
\end{eqnarray*}
Hence, we find after straightforward calculations that
\begin{equation}
\label{M3EPS}
\dE[\varepsilon_{n+1}^3|\cF_n]=\Big( q+a\frac{S_n}{n} \Big)-3\Big( q+a\frac{S_n}{n} \Big)^2+2\Big( q+a\frac{S_n}{n} \Big)^3 
\hspace{1cm}\text{a.s.}
\notag
\end{equation}
and
\begin{equation}
\label{M4EPS}
\dE[\varepsilon_{n+1}^4|\cF_n]=\Big( q+a\frac{S_n}{n} \Big)-4\Big( q+a\frac{S_n}{n} \Big)^2+6\Big( q+a\frac{S_n}{n} \Big)^3 
-3\Big( q+a\frac{S_n}{n} \Big)^4
\hspace{1cm}\text{a.s.}
\notag
\end{equation}
which ensures that
\begin{equation}
\label{SUPEPS4}
\sup_{n\geq 0} \mathbb{E}\left[\varepsilon^4_{n+1}|\mathcal{F}_n\right]\leq \frac{1}{12} \hspace{1cm} \text{a.s.}
\end{equation}
Hereafter, we deduce from \eqref{IPM} and \eqref{M2EPS} that
\begin{equation}
\label{CALCIP}
\langle M \rangle_n =  (s-q)(1-2q)+q(1-q)v_n+a(1-2q)\xi_n
-a^2 \zeta_n
\hspace{1cm} \text{a.s.}
\end{equation}
where we have denoted
\begin{equation} \label{vn}
v_n=\sum_{k=1}^n a_k^2,
\end{equation}
$$
\xi_n= \sum_{k=1}^{n-1} a_{k+1}^2 \Big(\frac{S_k}{k} \Big)
\hspace{1cm} \text{and} \hspace{1cm} 
\zeta_n= \sum_{k=1}^{n-1} a_{k+1}^2 \Big(\frac{S_k}{k} \Big)^2.
$$
Using standard results on the asymptotic behavior of the Euler Gamma function, we obtain that the MRW has 
diffusive, critical and superdiffusive regimes, depending on whether 
$a<1/2$, $a=1/2$ and $a>1/2$, respectively. In the diffusive regime,
\begin{equation}
\label{vn-DR}
\lim_{n \rightarrow \infty} \frac{v_n}{n^{1-2a}}=\ell \hspace{1cm}\text{where}\hspace{1cm}
\ell=\frac{\Gamma^2(a+1)}{1-2a}.
\end{equation}
In the critical regime,
\begin{equation}
\label{vn-CR}
\lim_{n \rightarrow \infty} \frac{v_n}{\log n}=\frac{\pi}{4}.
\end{equation}
In the superdiffusive regime, $(v_n)$ converges to a finite value. More precisely, as in \cite{Bercu2018},
\begin{equation} 
\label{vn-SR}
\lim_{n \rightarrow \infty}  v_n =\sum_{k=0}^ \infty \Big(\frac{\Gamma(a+1)\Gamma(k+1)}{\Gamma(k+a+1)} \Big)^2 = 
{}_{3}F_2 \Bigl( \begin{matrix}
{\hspace{-0.2cm}1 \ , \ 1 \ ,\hspace{0.2cm}1}\\
{a+1,a+1}\end{matrix} \Bigl|
{\displaystyle 1}\Bigr)
\end{equation} 
where $\!{}_{p}F_q$ stands for the hypergeometric function defined for all $z \in \dC$ by
\begin{equation}
{}_{p}F_q \Bigl( \begin{matrix}
{a_1,\ldots,a_p}\\
{b_1,\ldots,b_q}\end{matrix} \Bigl|
{\displaystyle z}\Bigr)
=\sum_{n=0}^{\infty}
\frac{(a_1)_n\,\cdots\,(a_p)_n}
{(b_1)_n\,\cdots\,(b_q)_n\, n!} z^n.
\notag
\end{equation}
All the above convergences will be the keystones in order to investigate the asymptotic behavior for the MRW.

\vspace{-1ex}
\section*{Appendix B. The diffusive regime}
\renewcommand{\thesection}{\Alph{section}}
\renewcommand{\theequation}{\thesection.\arabic{equation}}
\setcounter{section}{2}
\setcounter{equation}{0}
\vspace{-1ex}

{\bfseries Proof of Theorem \ref{T-ASP-DR}}. It is only necessary to prove the quadratic strong law \eqref{QSL-DR}.
Denote by $f_n$ the explosion coefficient associated with the martingale $(M_n)$,
\begin{equation}
f_n=\frac{a_n^2}{v_n} \nonumber
\end{equation}
We clearly have from \eqref{an} and \eqref{vn-DR} that $f_n$ converges to zero almost surely as $n$ goes to infinity. 
Moreover, we deduce from \eqref{SLLN-DR} and \eqref{M2EPS} that
\begin{equation}
\label{CVGM2EPS}
\lim_{n\rightarrow \infty} \mathbb{E}\left[\varepsilon^2_{n+1}|\mathcal{F}_n\right] = \sigma^2
\hspace{1cm} \text{a.s.}
\end{equation}
where the asymptotic variance $\sigma^2$ is given by \eqref{VAR-DR}.
Convergence \eqref{CVGM2EPS} together with the upper bound \eqref{M4EPS} and Theorem $3$ in \cite{Bercu2004} lead us to
\begin{equation} 
\lim_{n\rightarrow \infty} \frac{1}{\log v_n}\sum_{k=1}^n f_k \Big(\frac{M_k^2}{v_{k}}\Big)=\sigma^2 \hspace{1cm}\text{a.s.}
\notag
\end{equation}
Hence, it follows from convergence \eqref{vn-DR}  that
\begin{equation}
\lim_{n\rightarrow \infty} \frac{1}{\log n}\sum_{k=1}^n \Big(\frac{a_k M_k}{v_k}\Big)^2=(1-2a)\sigma^2 \hspace{1cm}\text{a.s.}
\notag
\end{equation}
Additionally, we get from the definition of $(M_n)$ that
\begin{equation}
\lim_{n\rightarrow \infty} \frac{1}{\log n}\sum_{k=1}^n \frac{a_k^4}{v_k^2} \Big(S_k-\frac{ q A_k}{a_k}\Big)^2 =(1-2a)\sigma^2 \hspace{1cm}\text{a.s.}
\label{QSLP1}
\end{equation}
However, one can easily see from \eqref{an} and \eqref{vn-DR} that
\begin{equation}
\label{EQAV-DR}
\lim_{n\rightarrow \infty}  \frac{n^2 a_n^4}{v_n^2}=(1-2a)^2.
\end{equation}
Therefore, we deduce from \eqref{QSLP1} and \eqref{EQAV-DR} that
\begin{equation}
\label{QSLP2}
\lim_{n\rightarrow \infty} \frac{1}{\log n}\sum_{k=1}^n  \Big(\frac{S_k}{k}-\frac{q A_k}{k a_k}\Big)^2 
=\frac{\sigma^2}{1-2a} \hspace{1cm}\text{a.s.}
\end{equation}
Furthermore, we clearly have for all $n\geq 1$,
$$
\Big(\frac{S_n}{n}-\frac{q}{1-a}\Big)^2\!\!=\!\Big(\frac{S_n}{n}-\frac{q A_n}{n a_n}\Big)^2\!\!+\!\Big(\frac{q A_n}{n a_n}-\frac{q}{1-a}\Big)^2
\!\!+\! 2\Big(\frac{S_n}{n}-\frac{q A_n}{n a_n}\Big)\!\Big(\frac{q A_n}{n a_n}-\frac{q}{1-a}\Big).
$$
Hereafter, by virtue of Lemma B.$1$ in \cite{Bercu2018}, we obtain that
\begin{equation}
\frac{A_n}{n a_n}=\frac{\Gamma(n+a)}{n\Gamma(n)}\sum_{k=1}^n \frac{\Gamma(k)}{\Gamma(k+a)}=
\frac{1}{a-1}\Big(\frac{\Gamma(n+a)}{\Gamma(n+1)\Gamma(a)}-1\Big)
\notag
\end{equation}
which implies that
\begin{equation}
\frac{A_n}{n a_n}-\frac{1}{1-a}=\frac{\Gamma(n+a)}{(a-1)\Gamma(a)\Gamma(n+1)}.
\label{EQUIFRACAN}
\end{equation}
Consequently, we find from \eqref{EQUIFRACAN}  that
\begin{equation}
\label{QSLP3}
\lim_{n\rightarrow \infty} \sum_{k=1}^n  \Big(\frac{A_k}{k a_k}-\frac{1}{1-a}\Big)^2
=\frac{1}{(a-1)^2\Gamma(a)^2} \sum_{k=1}^\infty \frac{1}{k^{2(1-a)}}
\end{equation}
which is finite because $2(1-a)>1$. Thus, \eqref{QSLP3} yields to
\begin{equation}
\label{QSLP4}
\lim_{n\rightarrow \infty} \frac{1}{\log n}\sum_{k=1}^n  \Big(\frac{A_k}{k a_k}-\frac{1}{1-a}\Big)^2=0.
\end{equation}
Finally, it follows from the Cauchy-Schwarz inequality together with
the almost sure convergence \eqref{QSLP2} and \eqref{QSLP4} that 
\begin{equation}
\lim_{n\rightarrow \infty} \frac{1}{\log n}\sum_{k=1}^n \Big(\frac{S_k}{k}-\frac{q}{1-a} \Big)^2=\frac{\sigma^2}{1-2a} \hspace{1cm}\text{a.s.}
\notag
\end{equation}
which completes the proof of Theorem \ref{T-ASP-DR}.
\demend


\noindent
{\bfseries Proof of Theorem \ref{T-FCLT-DR}}. We shall now proceed to the proof of the functional central limit theorem 
given by the distributional convergence \eqref{FCLT-DR}. It follows from \eqref{SLLN-DR}, \eqref{CALCIP}, \eqref{vn-DR}, together
with Toeplitz lemma \cite{Duflo1997}, that
\begin{equation}
\lim_{n\rightarrow \infty} \frac{1}{n^{1-2a}}\langle M \rangle_n =\sigma^2\ell \hspace{1cm}\text{a.s.}
\label{CVGIP-DR}
\end{equation}
where the asymptotic variance $\sigma^2$ is given by \eqref{VAR-DR}. Consequently, we deduce from \eqref{CVGIP-DR} that for all $t \geq 0$, 
\begin{equation}
\lim_{n\rightarrow \infty} \frac{1}{n^{1-2a}}\langle M \rangle_{\lfloor nt \rfloor} =\sigma^2\ell t^{1-2a} \hspace{1cm}\text{a.s.}
\label{CVGIPT-DR}
\end{equation}
It is now necessary to check that Lindeberg's condition is satisfied. In other words, we have to prove that for any $\eta>0$,
\begin{equation}
\label{LINDEBERG-DR}
\frac{1}{n^{1-2a}}\sum_{k=1}^{n}\dE\big[\Delta M_k^2 \rI_{\{|\Delta M_k|>\eta \sqrt{n^{1-2a}} \}}|\cF_{k-1}\big] \limp 0
\end{equation}
where $\Delta M_n = M_n - M_{n-1}$. We obtain from bound \eqref{M4EPS} that for any $\eta>0$,
\begin{eqnarray}
\hspace{-1.4cm} \frac{1}{n^{1-2a}}\sum_{k=1}^{n}\dE\big[\Delta M_k^2 \rI_{\{|\Delta M_k|>\eta \sqrt{n^{1-2a}} \}}|\cF_{k-1}\big] 
&\leq& 
\frac{1}{n^{2(1-2a)} \eta^2} \sum_{k=1}^{n} \dE\big[\Delta M_k^4|\cF_{k-1}\big], \notag\\
&\leq&
\frac{1}{n^{2(1-2a)} \eta^2} \sum_{k=1}^{n} a_k^4\dE\big[\veps_k^4|\cF_{k-1}\big], \notag\\
&\leq&
\frac{1}{12 n^{2(1-2a)} \eta^2}\sum_{k=1}^{n} a_k^4.
\label{LIND-DR}
\end{eqnarray}
However, it follows from convergence \eqref{vn-DR} together with \eqref{EQAV-DR} that
\begin{equation}
\lim_{n\rightarrow \infty}  \frac{1}{n^{1-2a}} \sum_{k=1}^{n} a_k^4=(1-2a) \ell^2.
\notag
\end{equation}
Hence, \eqref{LIND-DR} ensures that Lindeberg's condition is satisfied. Therefore, we immediately deduce from \eqref{LINDEBERG-DR} that
for all $t \geq 0$ and for any $\eta>0$,
\begin{equation}
\label{LINDEBERGT-DR}
\frac{1}{n^{1-2a}}\sum_{k=1}^{\lfloor nt \rfloor}\dE\big[\Delta M_k^2 \rI_{\{|\Delta M_k|>\eta \sqrt{n^{1-2a}} \}}|\cF_{k-1}\big] \limp 0.
\end{equation}
Consequently, we obtain from the functional central limit theorem  for martingales given in Theorem 2.5 of \cite{Durrett1978} that
\begin{equation}
\label{FCLTMART-DR}
\Big(\frac{M_{\lfloor nt \rfloor}}{\sqrt{n^{1-2a}}}, t \geq 0\Big) \Longrightarrow \big( B_t, t \geq 0 \big)
\end{equation}
where $\big( B_t, t \geq 0 \big)$ is a real-valued centered Gaussian process starting at the origin with covariance given, for all $0<s \leq t$, by
$\dE[B_s B_t]= \sigma^2 \ell s^{1-2a}$.
Hereafter, we find from the definition of $(M_n)$ together with \eqref{an} and \eqref{EQUIFRACAN} that
\begin{equation}
\label{DECFCLT-DR}
\frac{M_{\lfloor nt \rfloor}}{\sqrt{n^{1-2a}}}=\frac{\lfloor nt \rfloor a_{\lfloor nt \rfloor}}{\sqrt{n^{1-2a}}}\Big( \frac{S_{\lfloor nt \rfloor}}{\lfloor nt \rfloor}
- \frac{q}{1-a}\Big) - \frac{aq}{(a-1) \sqrt{n^{1-2a}}}.
\end{equation}
The right-hand side of \eqref{DECFCLT-DR} clearly goes to zero as $n$ tends to infinity. Furthermore, we obtain once again from \eqref{an} that
$$
\lim_{n\rightarrow \infty} \frac{\lfloor nt \rfloor a_{\lfloor nt \rfloor}}{n^{1-a}}= t^{1-a} \Gamma(a+1)
$$
Finally, we deduce from \eqref{DECFCLT-DR} that
\begin{equation}
\left( \sqrt{n}\Big(\frac{S_{\lfloor nt \rfloor}}{\lfloor nt \rfloor}-\frac{q}{1-a}\Big), t \geq 0\right) \Longrightarrow \big( W_t, t \geq 0 \big)
\notag
\end{equation}
where 
$$
W_t=\frac{B_t}{t^{1-a}\Gamma(a+1)},
$$
which means that
$\big( W_t, t \geq 0 \big)$ is a real-valued centered Gaussian process starting at the origin with covariance given, for all $0<s \leq t$, by
$$
\dE[W_s W_t]= \frac{\dE[B_s B_t]}{(st)^{1-a}\Gamma^2(a+1)}=\Big(\frac{\sigma^2}{1-2a} \Big) \frac{s^{1-2a}}{(st)^{1-a}}=\frac{\sigma^2}{(1-2a)t} \Big(\frac{t}{s}\Bigr)^a.
$$
\demend

\noindent
{\bfseries Proof of Corollary \ref{C-CM-DR}}. The almost sure convergence \eqref{SLLN-CM-DR} immediately follows from
\eqref{SLLN-DR} together with Toeplitz lemma \cite{Duflo1997}. Moreover, one can observe that the center of mass $G_n$ satisfies
$$
G_n=\int_0^1 S_{{\lfloor nt \rfloor}}dt.
$$
Consequently, the random variable
$$
\sqrt{n}\Big(\frac{G_n}{n}-\frac{q}{2(1-a)}\Big) 
$$
shares the same asymptotic normality than the continuous functional
$$
\int_0^1 \sqrt{n}\Big(\frac{S_{\lfloor nt \rfloor}}{\lfloor nt \rfloor}-\frac{q}{1-a}\Big) t dt.
$$
Applying Theorem \ref{T-FCLT-DR}, we immediately obtain that
\begin{equation}
\label{CLTGNP1}
\sqrt{n}\Big(\frac{G_n}{n}-\frac{q}{2(1-a)}\Big) \liml \int_0^1 t W_t dt.
\end{equation}
The right-hand side of \eqref{CLTGNP1} is a Gaussian random variable with zero mean and variance given by
\begin{eqnarray*}
\hspace{-1cm}
\dE\Big[\Big( \int_0^1 t W_t dt\Big)^2\Big]&=&2 \int_0^1 \int_0^t st \dE[W_s W_t] ds dt=\frac{2 \sigma^2}{1-2a} \int_0^1 \int_0^t s \Big(\frac{t}{s}\Big)^a ds dt,\\
&=& \frac{2 \sigma^2}{1-2a} \int_0^1 t^a \Big(\int_0^t s^{1-a} ds \Big)dt= \frac{2 \sigma^2}{(1-2a)(2-a)} \int_0^1 t^2 dt, \\
&=& \frac{2 \sigma^2}{3(1-2a)(2-a)},
\end{eqnarray*}
which completes the proof of Corollary \ref{C-CM-DR}.
\demend

\vspace{-1ex}
\section*{Appendix C. The critical regime}
\renewcommand{\thesection}{\Alph{section}}
\renewcommand{\theequation}{\thesection.\arabic{equation}}
\setcounter{section}{3}
\setcounter{equation}{0}
\vspace{-1ex}

{ \bfseries Proof of Theorem \ref{T-ASP-CR}}. 
It is only necessary to prove the quadratic strong law \eqref{QSL-CR}.
We have from \eqref{an} and \eqref{vn-CR} that $f_n$ converges to zero almost surely as $n$ goes to infinity. 
Moreover, we obtain from \eqref{SLLN-CR} and \eqref{M2EPS} that
\begin{equation}
\label{CVGM2EPS-CR}
\lim_{n\rightarrow \infty} \mathbb{E}\left[\varepsilon^2_{n+1}|\mathcal{F}_n\right] = 4q(1-p)
\hspace{1cm} \text{a.s.}
\end{equation}
Hence, it follows from Theorem $3$ in \cite{Bercu2004} and \eqref{CVGM2EPS-CR} that
\begin{equation} 
\lim_{n\rightarrow \infty} \frac{1}{\log v_n}\sum_{k=1}^n f_k \Big(\frac{M_k^2}{v_{k}}\Big)=4q(1-p) \hspace{1cm}\text{a.s.}
\notag
\end{equation}
which implies via \eqref{vn-CR} that
\begin{equation}
\lim_{n\rightarrow \infty} \frac{1}{\log \log n}\sum_{k=1}^n \Big(\frac{a_k M_k}{v_k}\Big)^2=4q(1-p) \hspace{1cm}\text{a.s.}
\label{QSLP5}
\end{equation}
However, one can observe that from \eqref{an} and \eqref{vn-CR} that
\begin{equation}
\label{EQAV-CR}
\lim_{n\rightarrow \infty}  \frac{(n \log n)^2 a_n^4}{v_n^2}=1.
\end{equation}
Consequently, we deduce from  \eqref{QSLP5} and \eqref{EQAV-CR} together with the definition of $(M_n)$ that
\begin{equation}  
\label{QSLP6}
\lim_{n\rightarrow \infty} \frac{1}{\log \log n}\sum_{k=2}^n \Big(\frac{1}{\log k}\Big)^2 \Big(\frac{S_k}{k}-\frac{qA_k}{ka_k} \Big)^2=4q(1-p) \hspace{1cm}\text{a.s.}
\end{equation}
Furthermore, we obtain from \eqref{EQUIFRACAN} that
\begin{equation}
\lim_{n\rightarrow \infty} \sum_{k=2}^n \Big(\frac{1}{\log k}\Big)^2 \Big(\frac{A_k}{k a_k}-2\Big)^2 = \frac{4}{\pi}\sum_{k=2}^\infty \frac{1}{k (\log k)^2} < + \infty.
\notag
\end{equation}
It ensures that 
\begin{equation}
\lim_{n\rightarrow \infty} \frac{1}{\log \log n} \sum_{k=2}^n \Big(\frac{1}{\log k}\Big)^2 \Big(\frac{A_k}{k a_k}-2\Big)^2 =0
\label{QSLP7}
\end{equation}
Finally, we find from the Cauchy-Schwarz inequality together with the almost sure convergence \eqref{QSLP6} and \eqref{QSLP7} that 
\begin{equation}
\lim_{n\rightarrow \infty} \frac{1}{\log \log n}\sum_{k=2}^n \Big(\frac{1}{\log k}\Big)^2 \Big(\frac{S_k}{k}-2q \Big)^2= 4q(1-p) \hspace{1cm}\text{a.s.}
\notag
\end{equation}
which is exactly what we wanted to prove.
\demend

\noindent
\hspace{-1ex}
{ \bfseries Proof of Theorem \ref{T-FCLT-CR}}. The proof of the distributional convergence 
\eqref{FCLT-CR} is left to the reader as it follows essentially the same lines as that of \eqref{FCLT-DR}.
The only minor change is that the almost sure rate of convergence $n^{1-2a}$ has to be replaced by 
$\log n$.
\demend


\noindent
{\bfseries Proof of Corollary \ref{C-CM-CR}}. The almost sure convergence \eqref{SLLN-CM-CR} immediately follows from \eqref{SLLN-CR} together with Toeplitz lemma \cite{Duflo1997}. The proof of 
the asymptotic normality \eqref{CLT-CM-CR} is left to the reader as it is similar to that of Theorem 2.6 in \cite{Bercu2021}.
\demend
\vspace{-1ex}
\section*{Appendix D. The superdiffusive regime}
\renewcommand{\thesection}{\Alph{section}}
\renewcommand{\theequation}{\thesection.\arabic{equation}}
\setcounter{section}{4}
\setcounter{equation}{0}
\vspace{-1ex}

\noindent
{ \bfseries Proof of Theorem \ref{T-ASP-SR}}. We now proceed to the proof of Theorem \ref{T-ASP-SR}.
For that purpose, we claim that the martingale $(M_n)$ is bounded in $\mathbb{L}^2$. As a matter of
fact, we have from the martingale decomposition \eqref{martingale2} together with
\eqref{SUPEPS2} that for all $n\geq 1$, 
\begin{equation}
\dE[M_{n+1}^2|\cF_n]=M_n^2+a_{n+1}^2\dE[\varepsilon_{n+1}^2|\cF_n] \leq M_n^2+\frac{a_{n+1}^2}{4} 
\hspace{1cm}\text{a.s.}
\notag
\end{equation}
leading to
$$\dE[M_{n+1}^2]\leq\dE[M_{n}^2]+\frac{a_{n+1}^2}{4}.$$
Consequently, it follows from \eqref{vn-SR} that
$$\sup_{n\geq 1} \dE[M_{n}^2]\leq \dE[M_1^2]  + \frac{v_n}{4}
\leq (1-q)^2 + \frac{1}{4}
{}_{3}F_2 \Bigl( \begin{matrix}
{\hspace{-0.2cm}1 \ , \ 1 \ ,\hspace{0.2cm}1}\\
{a+1,a+1}\end{matrix} \Bigl|
{\displaystyle 1}\Bigr),
$$
which means that the martingale $(M_n)$ is bounded in $\mathbb{L}^2$. Hence,
$(M_n)$ converges almost surely and in $\mathbb{L}^2$ to the random variable
\begin{equation}
\label{DEFM}
M=\sum_{k=1}^\infty a_k \epsilon_k.
\end{equation}
Hereafter, as $M_n=a_n S_n-q A_n$, we clearly find that
\begin{equation}
\label{ASCVGSRP}
\lim_{n\rightarrow \infty} n a_n \Big( \frac{S_n}{n} - \frac{qA_n}{n a_n} \Big)= M 
\hspace{1cm}\text{a.s.}
\end{equation}
However, we already saw from \eqref{EQUIFRACAN} that
\begin{equation}
\label{EXACTAN}
\frac{A_n}{na_n}=\frac{1}{1-a}-\frac{a}{(1-a)na_n}.
\end{equation}
Therefore, we deduce from \eqref{an}, \eqref{ASCVGSRP} and \eqref{EXACTAN} that
\begin{equation}
\label{SLLN-SRP}
\lim_{n\rightarrow \infty}
n^{1-a}\Big( \frac{S_n}{n}-\frac{q}{1-a}\Big) = L \hspace{1cm}\text{a.s.}
\end{equation}
where the limiting random variable $L$ is given by
\begin{equation}
\label{DEFL}
L=\frac{1}{\Gamma(a+1)}\Big(M-\frac{qa}{1-a}\Big).
\end{equation}
One can observe that the almost sure convergence \eqref{SLLN-SRP} clearly implies \eqref{FSLLN-SR}.
It only remains to prove the mean square convergence \eqref{SDL2-SR}. Since 
$(M_n)$ converges to $M$  in $\mathbb{L}^2$, we have
${\displaystyle\lim_{n\rightarrow \infty}} \dE[(a_n S_n -q A_n -M)^2]=0$, which implies that
\begin{equation}
\label{L2SP}
\lim_{n\rightarrow \infty} \dE\Big[\Big( na_n\Big(\frac{S_n}{n}-\frac{qA_n}{n a_n} -M \Big)^{\! 2}\Big]=0.
\end{equation}
Dividing on both sides by $\Gamma^2(a+1)$, it follows from
\eqref{L2SP} together with  \eqref{EXACTAN} and \eqref{DEFL} that
$$\lim_{n\rightarrow \infty} 
\dE\Big[\Big(n^{1-a}\Big(\frac{S_n}{n}-\frac{q}{1-a}\Big) -L \Big)^2 \Big]=0,
$$
which leads us to
$$\lim_{n\rightarrow \infty} 
\dE\Big[\Big(\frac{n a_n}{\Gamma(a+1)}\Big(\frac{S_n}{n}-\frac{q}{1-a}\Big) -L \Big)^2 \Big]=0,
$$
which completes the proof of Theorem \ref{T-ASP-SR}.
\demend

\noindent
{ \bfseries Proof of Theorem \ref{T-MOML}}. 
By taking expectations on both sides of \eqref{expos}, we find that for all $n\geq 1$, 
\begin{equation}
\mathbb{E}[S_{n+1}]=q+\Big(1+\frac{a}{n}\Big)\mathbb{E}[S_n].
\notag
\end{equation}
Consequently, 
\begin{equation}
\label{MSN1}
\mathbb{E}[S_n]=\mathbb{E}[X_1]\prod_{k=1}^n \Big(1+\frac{a}{k}\Big)+q\sum_{k=1}^{n-1} \prod_{i=k+1}^{n-1} \Big(1+\frac{a}{i}\Big).
\end{equation}
Hence, we obtain from \eqref{an}, \eqref{bigAn} and \eqref{MSN1} that
\begin{equation}
\label{MSN2}
a_n\dE[S_n]=s-q+qA_n.
\end{equation}
Therefore, \eqref{MSN2} ensures that
\begin{equation}
\lim_{n\rightarrow \infty} \dE[M_n]=\dE[M]=s-q,
\label{MM1}
\end{equation}
leading to
$$
\dE[L]=\frac{1}{\Gamma(a+1)}\Big(\dE[M]-\frac{qa}{1-a}\Big)
=
\frac{s+\nu}{\Gamma(a+1)}
$$
where $\nu=-q/(1-a)$. We now compute the second moment of $L$.
By taking expectations on both sides of \eqref{Snsqq}, we obtain that for all $n\geq 1$, 
\begin{equation}
\dE[S_{n+1}^2]=\Big(1+\frac{2a}{n}\Big)\dE[S_n^2]+\Big(2q+\frac{a}{n}\Big)\dE[S_n]+q. 
\label{MSN3}
\end{equation}
If we set
\begin{equation}
g_n=\Big(1+\frac{2a}{n}\Big) \hspace{1cm}\text{and}\hspace{1cm} h_n=\Bigl(2q+\frac{a}{n}\Big)\dE[S_n]+q, 
\end{equation}
we have from \eqref{MSN3} that for all $n\geq 1$,
$\mathbb{E}[S_{n+1}^2]=g_n\mathbb{E}[S_{n}^2]+h_n$, which implies that
\begin{equation} 
\label{MSN4}
\mathbb{E}[S_{n}^2]=\frac{\Gamma(n+2a)}{\Gamma(n)\Gamma(2a+1)}\Big(s+\sum_{k=1}^{n-1}h_k\frac{\Gamma(2a+1)\Gamma(k+1)}{\Gamma(k+1+2a)}\Big).
\end{equation}
However, it is not hard to see from \eqref{EXACTAN} and \eqref{MSN2} that
\begin{eqnarray*}
h_n &=&  q+\Bigl(2q+\frac{a}{n}\Big) \Big(\frac{s-q+qA_n}{a_n}\Big), \\
&=& q+\frac{(2qn+a)(s-q)}{n a_n} + \frac{q(2qn+a)A_n}{na_n}, \\
&=& \frac{q}{1-a}+\frac{2q^2 n}{1-a}-\frac{b(2qn+a)}{(1-a)n a_n},
\end{eqnarray*}
where $b=q-s(1-a)$. Consequently, the strategy for finding a more tractable expression for $\mathbb{E}[S_n^2]$ deals with the simplification of
\begin{eqnarray*}
\hspace{-2cm}
\sum_{k=1}^{n-1}\frac{h_k\Gamma(k+1)}{\Gamma(k+1+2a)}
&=&\frac{q}{1-a}\sum_{k=1}^{n-1}\frac{\Gamma(k+1)}{\Gamma(k+1+2a)}+\frac{2q^2}{1-a}\sum_{k=1}^{n-1}\frac{k\Gamma(k+1)}{\Gamma(k+1+2a)}\\
&-&\frac{2bq}{1-a}\sum_{k=1}^{n-1}\frac{\Gamma(k+1)}{a_k\Gamma(k+1+2a)}-\frac{ab}{1-a}\sum_{k=1}^{n-1}\frac{\Gamma(k+1)}{k a_k\Gamma(k+1+2a)}.
\end{eqnarray*}
We shall now make repeatedly use of Lemma B.$1$ in \cite{Bercu2018} jointly with \eqref{an} and standard properties of the Gamma function. First of all,
\begin{eqnarray}
\hspace{-1cm}
\sum_{k=1}^{n-1}\frac{\Gamma(k+1)}{\Gamma(k+1+2a)}&=&\sum_{k=1}^{n}\frac{\Gamma(k)}{\Gamma(k+2a)}-\frac{1}{\Gamma(1+2a)} \notag\\
&=&\frac{1}{(2a-1)\Gamma(1+2a)}-\frac{n\Gamma(n)}{(2a-1)\Gamma(n+2a)}. \label{hn1}
\end{eqnarray}
Next, we also have
\begin{eqnarray}
\hspace{-2cm}
\sum_{k=1}^{n-1}\frac{k\Gamma(k+1)}{\Gamma(k+1+2a)}&=&\sum_{k=2}^{n}\frac{(k-1)\Gamma(k)}{\Gamma(k+2a)}=\sum_{k=1}^{n}\frac{\Gamma(k+1)}{\Gamma(k+2a)}-\sum_{k=1}^{n}\frac{\Gamma(k)}{\Gamma(k+2a)}\notag \\
&=&\frac{1}{2(a-1)(2a-1)\Gamma(2a)}-\frac{n\Gamma(n)\big(n(2a-1)+1\big)}{2(a-1)(2a-1)\Gamma(n+2a)}. \label{hn2}
\end{eqnarray}
In the same vein,
 \begin{eqnarray}
& \sum_{k=1}^{n-1}&\frac{\Gamma(k+1)}{a_k\Gamma(k+1+2a)}=\sum_{k=1}^{n-1}\frac{\Gamma(k+1)}{\Gamma(k+1+2a)}\frac{\Gamma(k+a)}{\Gamma(k)\Gamma(a+1)}\notag \\
&=&\frac{1}{\Gamma(a+1)}\sum_{k=1}^{n}\frac{(k-1)\Gamma(k-1+a)}{\Gamma(k+2a)} \notag \\
&=&\frac{1}{\Gamma(a+1)}\Big( \frac{\Gamma(a+1)}{(a-1)\Gamma(2a)}-\frac{\Gamma(n+a+1)}{(a-1)\Gamma(n+2a)}-\frac{\Gamma(a)}{\Gamma(2a)}+\frac{\Gamma(n+a)}{\Gamma(n+2a)} \Big)  \notag \\
&=&\frac{1}{a(a-1)\Gamma(2a)}-\frac{(n+1)\Gamma(n+a)}{a (a-1)\Gamma(a)\Gamma(n+2a)}. \label{hn3}
\end{eqnarray}
Furthermore,
\begin{eqnarray}
&\sum_{k=1}^{n-1}& \frac{\Gamma(k+1)}{k a_k\Gamma(k+1+2a)}=\sum_{k=1}^{n-1}\frac{k\Gamma(k)\Gamma(k+a)}{k \Gamma(k) \Gamma(a+1) \Gamma(k+1+2a)} \notag \\
&=&\frac{1}{\Gamma(a+1)}\sum_{k=2}^n \frac{\Gamma(k-1+a)}{\Gamma(k+2a)}\notag\\
&=&\frac{1}{\Gamma(a+1)} \Big(\frac{\Gamma(a)}{a\Gamma(2a)}-\frac{\Gamma(a)}{2a\Gamma(2a)}-\frac{\Gamma(n+a)}{a\Gamma(n+2a)}\Big)\notag \\
&=&\frac{1}{a\Gamma(2a+1)}-\frac{\Gamma(n+a)}{a^2\Gamma(a)\Gamma(n+2a)}. \label{hn4}
\end{eqnarray}
Putting together the four contributions \eqref{hn1} to \eqref{hn4}, we obtain after some tedious but 
straightforward calculations that
%

\begin{eqnarray}
\hspace{-1cm}
\dE[S_n^2]&=&\frac{(s+\tau)\Gamma(n+2a)}{\Gamma(n)\Gamma(2a+1)}
-\frac{qn}{(1-a)(2a-1)}+\frac{q^2n\big(n(2a-1) -1\big)}{(1-a)^2 (2a-1)} \notag\\
&+&\frac{b\Gamma(n+a) \big(-2q(n+1)+1-a \big)}{(1-a)^2\Gamma(n)\Gamma(a+1)}
\notag
\end{eqnarray}
where $\tau$ is given by \eqref{deftau}, leading to
\begin{eqnarray}
\hspace{-1cm}
a_n^2\dE[S_n^2]&=&\frac{(s+\tau)a_n^2\Gamma(n+2a)}{\Gamma(n)\Gamma(2a+1)}
-\frac{qn a_n^2}{(1-a)(2a-1)}+\frac{q^2n a_n^2\big(n(2a-1) -1\big)}{(1-a)^2 (2a-1)} \notag\\
&+&\frac{b a_n \big(-2q(n+1)+1-a \big)}{(1-a)^2}.
\label{MSN5}
\end{eqnarray}
Herefater, it follows from the definition of $M_n$ together with \eqref{MSN2} that
$$
\dE[M_n^2]=a_n^2\dE[S_n^2]-2qa_n A_n\dE[S_n]+q^2A_n^2=
a_n^2\dE[S_n^2]-2q(s-q)A_n -q^2A_n^2.
$$
Consequently, we deduce from \eqref{EXACTAN} and \eqref{MSN5} that
\begin{eqnarray}
\hspace{-1cm}
\dE[M_n^2]&=&\frac{(s+\tau)a_n^2\Gamma(n+2a)}{\Gamma(n)\Gamma(2a+1)}
-\frac{qn a_n^2}{(1-a)(2a-1)}+\frac{q^2n^2 a_n^2}{(1-a)^2} \notag\\
&-&\frac{q^2n a_n^2}{(1-a)^2(2a-1)} + \frac{b a_n(1-a-2q)}{(1-a)^2} -\frac{2bqna_n}{(1-a)^2} -\frac{2q(s-q)n a_n}{1-a} \notag \\
&+& \frac{2qa(s-q)}{1-a} 
-\frac{q^2 n^2 a_n^2}{(1-a)^2}-\frac{q^2a^2}{(1-a)^2}+\frac{2q^2 a n a_n}{(a-1)^2}.
\label{MSN6}
\end{eqnarray}
Therefore, as $b+(s-q)(1-a) -aq=0$, \eqref{MSN6} reduces to
\begin{eqnarray}
\dE[M_n^2]&=&\frac{(s+\tau)a_n^2\Gamma(n+2a)}{\Gamma(n)\Gamma(2a+1)}
-\frac{qa\big( q(2-a)-2s(1-a) \big)}{(1-a)^2} \notag\\
&-&\frac{qn a_n^2(1-a+q)}{(1-a)(2a-1)}+ \frac{b a_n(1-a-2q)}{(1-a)^2}.
\label{MSN7}
\end{eqnarray}
Hence, we deduce from \eqref{an} and \eqref{MSN7} that
\begin{equation}
\lim_{n\rightarrow \infty} \dE[M_n^2]=\dE[M^2]=\frac{(s+\tau)\Gamma^2(a+1)}{\Gamma(2a+1)}
-\frac{qa\big( q(2-a)-2s(1-a) \big)}{(1-a)^2}
\label{MM2}
\end{equation}
Finally, we find from \eqref{MM1} and \eqref{MM2} that
\begin{equation}
\dE[L^2]=\frac{1}{\Gamma^2(a+1)}\Big(\dE[M^2]-\frac{2qa}{1-a}\dE[M]+\frac{q^2 a^2}{(1-a)^2}\Big)
=\frac{s+\tau}{\Gamma(2a+1)},
\notag
\end{equation}
which is exactly what we wanted to prove.
\demend

\noindent
{ \bfseries Proof of Theorem \ref{T-AN-SR}}. It only remains to prove the asymptotic normality
\eqref{AN-SR}. On the one hand, 
we clearly have from \eqref{SUPEPS2} and \eqref{vn-SR} that
$$
\sum_{n=1}^\infty \dE[\Delta M_n^2] \leq \dE[M_1^2]  + \frac{v_n}{4}
\leq (1-q)^2 + \frac{1}{4}
{}_{3}F_2 \Bigl( \begin{matrix}
{\hspace{-0.2cm}1 \ , \ 1 \ ,\hspace{0.2cm}1}\\
{a+1,a+1}\end{matrix} \Bigl|
{\displaystyle 1}\Bigr) < +\infty.
$$
On the other hand, let
$$
s_n^2=\sum_{k=n}^\infty \dE[\Delta M_k^2]=\sum_{k=n}^\infty a_k^2 \dE[\varepsilon_k^2].
$$
We claim that
\begin{equation}
\label{AN-SRP1}
\lim_{n\rightarrow \infty} n^{2a-1} s_n^2= \frac{q(1-q-a)}{2a-1}\Gamma^2(a+1)= \frac{\sigma^2 (1-a)^2}{2a-1}\Gamma^2(a+1)
\end{equation}
where the asymptotic variance $\sigma^2$ is given by \eqref{VAR-DR}.
As a matter of fact, we already saw from \eqref{M2EPS} that for all $n \geq 1$,
$$
\mathbb{E}[\varepsilon^2_{n+1}|\mathcal{F}_n]=
\Big(q+a\frac{S_n}{n} \Big)-\Big(q+a\frac{S_n}{n} \Big)^2 \hspace{1cm} \text{a.s.}
$$
which implies that for all $n \geq 1$,
\begin{equation}
\label{AN-SRP2}
\mathbb{E}[\varepsilon^2_{n+1}]=q(1-q)+(1-2q) \frac{a\dE[S_n]}{n} - \frac{a^2\dE[S_n^2]}{n^2}.
\end{equation}
Consequently, \eqref{AN-SRP1} follows from \eqref{MSN2}, \eqref{MSN5} and \eqref{AN-SRP2} 
together with the four convergences that arise directly from \eqref{an},
\begin{eqnarray*}
\lim_{n\rightarrow \infty} n^{2a-1} \sum_{k=n}^\infty a_k^2&=&\frac{\Gamma^2(a+1)}{2a-1},
\hspace{1.5cm} \lim_{n\rightarrow \infty} n^{a} \sum_{k=n}^\infty \frac{a_k}{k} =\frac{\Gamma(a+1)}{a}, \\
\lim_{n\rightarrow \infty} n^{2a} \sum_{k=n}^\infty \frac{a_k^2}{k}&=&\frac{\Gamma^2(a+1)}{2a},
\hspace{1.5cm} \lim_{n\rightarrow \infty} n^{a+1} \sum_{k=n}^\infty \frac{a_k}{k^2} =
\frac{\Gamma(a+1)}{a+1}.
\end{eqnarray*}
In addition, we deduce from \eqref{CVGM2EPS} and \eqref{AN-SRP1} that
\begin{equation}
\label{AN-SRP3}
\lim_{n\rightarrow \infty} \frac{1}{s_n^2} \sum_{k=n}^\infty  \dE[\Delta M_k^2| \cF_{k-1}]
=\lim_{n\rightarrow \infty} \frac{1}{s_n^2} \sum_{k=n}^\infty  a_k^2 \dE[\veps_k^2| \cF_{k-1}]
= \frac{1}{(1-a)^2}
\hspace{1cm} \text{a.s.}
\end{equation}
Furthermore, we obtain from bound \eqref{M4EPS} that for any $\eta>0$,
\begin{eqnarray}
\hspace{-1.4cm} \frac{1}{s_n^2}\sum_{k=n}^{\infty}\dE\big[\Delta M_k^2 \rI_{\{|\Delta M_k|>\eta s_n \}}\big] 
&\leq& 
\frac{1}{s_n^{4} \eta^2} \sum_{k=n}^{\infty} \dE\big[\Delta M_k^4\big] \leq
\frac{1}{s_n^{4} \eta^2} \sum_{k=n}^{\infty} a_k^4\dE\big[\veps_k^4\big], \notag\\
&\leq&
\frac{1}{12 s_n^{4} \eta^2}\sum_{k=1}^{n} a_k^4.
\label{LIND-SR}
\end{eqnarray}
However, it is not hard to see from \eqref{an} that
\begin{equation}
\lim_{n\rightarrow \infty}  n^{4a-1} \sum_{k=n}^{\infty} a_k^4=\frac{\Gamma^4(a+1)}{4a-1}.
\label{EQS4AN}
\end{equation}
Hence, \eqref{AN-SRP1} together with \eqref{LIND-SR} and \eqref{EQS4AN}, ensure that for any
$\eta >0$,
\begin{equation}
\label{AN-SRP4}
\lim_{n\rightarrow \infty} \frac{1}{s_n^2} \sum_{k=n}^\infty  \dE\big[\Delta M_k^2 \rI_{\{|\Delta M_k|>\eta s_n \}}\big] 
= 0.
\end{equation}
All the conditions of Theorem 1 and Corollary 1 in \cite{Heyde1977} are satisfied, which leads to the asymptotic normality
\begin{equation}
\label{AN-SRP5}
\frac{M_n-M}{s_n} \liml \cN\Big(0,  \frac{1}{(1-a)^2} \Big).
\end{equation}
Hereafter, we obtain from the definition of $L$ given by \eqref{DEFL} together with \eqref{AN-SRP1} and \eqref{AN-SRP5}, that
\begin{equation}
\label{AN-SRP6}
\sqrt{\frac{(2a-1)n^{2a-1}}{\sigma^2 (1-a)^2}}\left( \frac{na_n}{\Gamma(a+1)}\Big( \frac{S_n}{n} - \frac{q}{1-a} \Big) -L \right) \liml
\cN\Big(0,  \frac{1}{(1-a)^2} \Big).
\end{equation}
Finally, as 
\begin{equation}
\label{AN-SRP7}
\frac{a_n}{\Gamma(a+1)}=\frac{1}{n^a}\Big(1+\frac{a(1-a)}{2n}+ O\Big(\frac{1}{n^2}\Big) \Big),
\end{equation}
we deduce from \eqref{AN-SRP6} and \eqref{AN-SRP7} that
\begin{equation}
\sqrt{n^{2a-1}}\left( n^{1-a}\Big(\frac{S_n}{n}-\frac{q}{1-a}\Big)-L \right) 
\liml \cN\Big(0,\frac{\sigma^2}{2a-1}\Big)
\notag
\end{equation}
which acheives the proof of Theorem \ref{T-AN-SR}. 
\demend

\noindent
{\bfseries Proof of Corollary \ref{C-CM-SR}}. It have from the definition of $G_n$ that
\begin{eqnarray}
nG_n & = & \sum_{k=1}^n S_k=\sum_{k=1}^n \Big(\frac{S_k}{k} -\frac{q}{1-a}\Big)k + \Big(\frac{q}{1-a} \Big) \sum_{k=1}^n k, \notag \\
& = &\sum_{k=1}^n k^{1-a}\Big(\frac{S_k}{k} -\frac{q}{1-a}\Big)k^a +\frac{qn(n+1)}{2(1-a)}. 
\label{CM-SRP1}
\end{eqnarray}
It follows from \eqref{SLLN-SR} together with Toeplitz lemma \cite{Duflo1997} that
\begin{equation}
\label{CM-SRP2}
\lim_{n\rightarrow \infty} \frac{1}{n^{1+a}} \sum_{k=1}^n k^{1-a}\Big(\frac{S_k}{k} -\frac{q}{1-a}\Big)k^a
=\frac{L}{1+a}
\hspace{1cm}\text{a.s.}
\end{equation}
Hence, we deduce from \eqref{CM-SRP1} and \eqref{CM-SRP2} that
$$
\lim_{n\rightarrow \infty} \frac{n}{n^{1+a}} \Big( G_n - \frac{q(n+1)}{2(1-a)} \Big)
=\frac{L}{1+a}
\hspace{1cm}\text{a.s.}
$$
which clearly leads to
$$
\lim_{n\rightarrow \infty} n^{1-a} \Big( \frac{G_n}{n} - \frac{q}{2(1-a)}\Big)
=\frac{L}{1+a}
\hspace{1cm}\text{a.s.}
$$ 
completing the proof of Corollary \ref{C-CM-SR}.  \demend


\end{document}